\documentclass[preprint,12pt]{elsarticle}

\usepackage{mathtools} 
\usepackage{url}
\usepackage[breaklinks]{hyperref}
\usepackage[headings]{fullpage}
\usepackage{datetime}  
\usepackage[english]{babel} 

\usepackage[full]{textcomp}
\usepackage[osf]{newtxtext} 
\usepackage{cabin} 
\usepackage[bigdelims,vvarbb]{newtxmath} 
\usepackage[cal=boondoxo]{mathalfa} 
\usepackage{zlmtt}
 
\usepackage{microtype}
 
\newcommand{\eps}{\varepsilon}

\newcommand{\dx}{\mathrm{d}} 
\newcommand{\A}{\mathcal{A}} 
\newcommand{\B}{\mathcal{B}} 
 
\newcommand{\I}{\mathcal{I}}

\newcommand{\N}{\mathbb{N}}
\newcommand{\R}{\mathbb{R}}

\newcommand{\M}{\mathfrak{M}}
\newcommand{\m}{\mathfrak{m}}
\newcommand{\gotht}{\mathfrak{t}}
\newcommand{\varemptyset}{\text{\O}}

\newcommand{\Odi}[1]{\Odip{}{#1}}
\newcommand{\Odig}[1]{\mathcal{O}\Bigl(#1\Bigr)}
   
\newcommand{\Odip}[2]{\mathcal{O}_{#1}\left(#2\right)}

\newcommand{\odip}[2]{{o}_{#1}\left(#2\right)}
\newcommand{\odi}[1]{\odip{}{#1}}
\newcommand{\odim}[1]{o\bigl(#1\bigr)}

\allowdisplaybreaks

\newtheorem{Theorem}{Theorem}
\newtheorem{Lemma}{Lemma} 
 
\newproof{pf}{Proof}

\begin{document}

\begin{frontmatter}
\title{A Diophantine approximation problem with two primes and one $k$-th power of a prime} 
 

\author[ag]{Alessandro Gambini\corref{cor1}}
\ead{a.gambini@unibo.it}
\cortext[cor1]{Corresponding author}
\address[ag]{Universit\`a di Parma, Dipartimento di Scienze Matematiche, 
  Fisiche e Informatiche, Parco Area delle Scienze 53/a, 43124 Parma, Italy.}
\author[al]{Alessandro Languasco}
\ead{alessandro.languasco@unipd.it }
\address[al]{Universit\`a di Padova,  Dipartimento di Matematica 
 ``Tullio Levi-Civita'', Via Trieste 63, 35121 Padova, Italy.}
\author[ag]{Alessandro Zaccagnini}
\ead{alessandro.zaccagnini@unipr.it}


\begin{abstract} 
We refine a result of the last two Authors of \cite{Languasco-Zaccagnini2016}
on a Diophantine approximation problem with two primes and a $k$-th power
of a prime which was only proved to hold for $1 < k < 4 / 3$.
We improve the $k$-range to $1<k\le 3$ by combining Harman's
technique on the minor arc with a suitable estimate for the $L^4$-norm of
the relevant exponential sum over primes.
\end{abstract} 

\end{frontmatter}

\section{Introduction}

This paper deals with an improvement of the result contained in
\cite{Languasco-Zaccagnini2016}, which is due to the last two Authors:
we refer to its introduction for a more thorough description of the
general Diophantine problem with prime variables.
Here we just recall that the goal is to prove that the inequality
\[
  \vert
    \lambda_1 p_1^{k_1} + \cdots + \lambda_r p_r^{k_r} - \omega
  \vert
  \le
  \eta,
\]
where $k_1$, \dots, $k_r$ are fixed positive numbers, $\lambda_1$,
\dots, $\lambda_r$ are fixed non-zero real numbers and $\eta > 0$ is
arbitrary, has infinitely many solutions in prime variables $p_1$,
\dots, $p_r$ for any given real number $\omega$, under as mild
Diophantine assumptions on $\lambda_1$, \dots, $\lambda_r$ as
possible.
In some cases, it is even possible to prove that the above inequality
holds when $\eta$ is a small negative power of the largest prime
occurring in it, usually when $1 / k_1 + \cdots + 1 / k_r$ is large
enough.

The problem tackled in \cite{Languasco-Zaccagnini2016} had $r = 3$,
$k_1 = k_2 = 1$, $k_3 = k \in (1, 4 / 3)$.
Assuming that $\lambda_1 / \lambda_2$ is irrational and that the
coefficients $\lambda_j$ are not all of the same sign, the last two
Authors proved that one can take
$\eta = \bigl( \max \{ p_1, p_2, p_3^k \} \bigr)^{-\phi(k) + \eps}$
for any fixed $\eps > 0$, where $\phi(k) = (4 - 3 k) / (10 k)$.
Our purpose in this paper is to improve on this result both in
the admissible range for $k$ and in the exponent,
replacing $\phi(k)$ by a larger value in the common range.
More precisely, we prove the following Theorem.

\begin{Theorem}
\label{main-thm}
Assume that $1 < k \le 3$, $\lambda_1$, $\lambda_2$ and $\lambda_3$ are
non-zero real numbers, not all of the same sign, that
$\lambda_1 / \lambda_2$ is irrational and let $\omega$ be 
a real number.
The inequality
\begin{equation}
\label{teorema}
  \vert
    \lambda_1 p_1 + \lambda_2 p_2 + \lambda_3 p_3^k - \omega
  \vert
  \le \bigl(\max \{ p_1, p_2, p_3^k \} \bigr)^{- \psi(k) + \eps}
\end{equation}
has infinitely many solutions in prime variables $p_1$, $p_2$, $p_3$
for any $\eps > 0$, where 
\begin{equation}
\label{psi-value}
  \psi(k)
  =
  \begin{cases}
    (3 - 2 k) / (6 k) & \text{if $1 < k \le \frac65$,} \\
    1 / 12            & \text{if $\frac65 < k \le 2$,} \\
    (3 -   k) / (6 k) & \text{if $2 < k < 3$,} \\
    1 / 24            & \text{if $k = 3$.}
  \end{cases}
\end{equation}
\end{Theorem}

We point out that in the common range $1 < k < 4 / 3$ we
have $\psi(k) > \phi(k)$.
We also remark that the strong bounds for the exponential sum $S_k$,
defined in \eqref{U_k} below, that recently became available for
integral $k$ (see Bourgain \cite{Bourgain2016} and Bourgain, Demeter
\& Guth \cite{BourgainDG2016}) are not useful in our problem.

The technique used to tackle this problem is the variant of the
circle method introduced in the 1940's by Davenport \& Heilbronn
\cite{davenport1946indefinite}, where the integration on a circle, or
equivalently on the interval $[0, 1]$, is replaced by integration on
the whole real line.
Our improvement is due to the use of the Harman technique on the minor
arc and to the fourth-power average for the exponential sum $S_k$ for
$k \ge 1$.

We thank the anonymous referee for an extremely careful reading of a
previous version of this paper.

\section{Outline of the proof}

Throughout this paper $p_i$ denotes a prime number, $k\ge1$ is a real number, $\eps$ 
is an arbitrarily small positive number whose value may vary depending on the 
occurrences and $\omega$ is a fixed real number. In order to prove that \eqref{teorema} 
has infinitely many solutions, it is sufficient to construct an increasing sequence $X_n$ 
that tends to infinity such that \eqref{teorema} has at least one solution with 
$\max \{ p_1, p_2, p_3^k \} \in[\delta X_n,X_n]$, with a fixed $\delta>0$
which depends
only on the choice of $\lambda_1$, $\lambda_2$ and $\lambda_3$.
Let $q$ be a denominator of a convergent to $\lambda_1/\lambda_2$ and
let $X_n=X$ (dropping the suffix $n$) run through the sequence
$X=q^3$.
The main quantities we will use are:
\begin{equation}
\label{U_k}
S_k(\alpha)=\!\!\!\sum_{\delta X\le p^k\le X} \!\!\!\log p\  e(p^k\alpha),\quad
U_k(\alpha)=\!\!\!\sum_{\delta X\le n^k\le X}\!\!\! e(n^k\alpha) \quad \textrm{and} \quad
T_k(\alpha)=\int_{(\delta X)^{1/k}}^{X^{1/k}}e( t^k\alpha)\,\dx  t,
\end{equation}
where $e(\alpha)=e^{2\pi i\alpha}$. We will approximate $S_k$ with $T_k$ and $U_k$.
By the Prime Number Theorem and   first derivative estimates for trigonometric integrals 
we have
\begin{equation}
\label{stima_tk}
S_k(\alpha)\ll_{k,\delta} X^{1/k}, \qquad
T_k(\alpha)\ll_{k,\delta} X^{1/k-1} \min\{ X, \vert \alpha\vert ^{-1}\}.
\end{equation}
Moreover the Euler summation formula implies that
\begin{equation}
\label{t-u}
T_k(\alpha)-U_k(\alpha)\ll_{k,\delta} 1+\vert \alpha\vert X.
\end{equation}

We also need a continuous function we will use to detect the solutions of \eqref{teorema},
 so we introduce
\[
  \widehat{K}_\eta(\alpha)
  :=
  \max\{ 0, \eta-\vert \alpha\vert \},
  \quad\textrm{where}\ \eta>0,
\]
which is the Fourier transform of the function $K_\eta$ defined by
\[
K_{\eta}(\alpha)=\Bigl(\frac{\sin(\pi\alpha\eta)}{\pi\alpha}\Bigr)^2
\]
for $\alpha\neq0$ and, by continuity, $K_{\eta}(0)=\eta^2$.
A well-known estimate is
\begin{equation}
\label{k_eta}
  K_{\eta}(\alpha) \ll \min\{ \eta^2, \vert \alpha\vert ^{-2} \}.
\end{equation}
 
Let now
\[
\mathcal{P}(X)
=
\{(p_1,p_2,p_3): \delta X<p_1,p_{2} \le X,\ \delta X<p_3^k \le X\}
\]
and
\[
\I(\eta,\omega,\mathfrak{X})
=
\int_{\mathfrak{X}}S_1(\lambda_1\alpha)S_1(\lambda_2\alpha)S_k(\lambda_3\alpha)
K_{\eta}(\alpha)e(-\omega\alpha)\, \dx \alpha,
\]
where $\mathfrak{X}$ is a measurable subset of $\R $.
{}From \eqref{U_k} and using the Fourier transform
 of $K_{\eta}(\alpha)$, we get
\begin{align*}
\notag
\I(\eta,\omega,\R )
&=
\sum_{(p_1, p_2, p_3) \in \mathcal{P}(X)}
\log p_1\log p_2\log p_3\ 
  \max\bigl\{0, \eta- \vert \lambda_1p_1+\lambda_2p_2+\lambda_3p_3^k-\omega \vert \bigr\}
\\ 
&
\leq \eta(\log X)^3\mathcal{N}(X),
\end{align*}
where $\mathcal{N}(X)$ actually denotes the number of solutions of the inequality 
\eqref{teorema} with $(p_1, p_2, p_3) \in \mathcal{P}(X)$. In other words 
$\I(\eta,\omega,\R )$ provides a lower bound for the quantity we are 
interested in; therefore it is sufficient  to prove 
that $\I(\eta,\omega,\R )>0$.

We now decompose $\R $ into subsets such that $\R =\M\,\cup\ \M^*\,\cup\ \m\,\cup\, \gotht$ 
where $\M$ is the major arc, $\M^*$ is the intermediate arc (which is
non-empty only for some values of $k$, see section
\ref{sect:inter-arc}), $\m$ is the minor arc and
$\gotht$ is the trivial arc. The decomposition is the following:
if $1 < k < 5/2$ we consider
\begin{align}
\notag
  \M
  &=
  [-P/X, P/X],
  &
  \M^*
  &=
  \varemptyset, \\ 
\label{def-arcs-k-basso}
  \m
  &=
  [P/X, R]\cup [-R,-P/X],
  &
  \gotht
  &=
  \R \setminus (\M\cup\M^*\cup \m),
\end{align}
while, for $5/2 \le k \le 3$, we set
\begin{align}
\notag
  \M
  &=
  \left[-P/X,P/X\right],
  &
  \M^*
  &=
  [P/X,X^{-3/5}] \cup [-X^{-3/5},-P/X ], \\
\label{def-arcs-k-alto}
  \m
  &=
  [X^{-3/5},R] \cup [-R,-X^{-3/5}],
  &
  \gotht
  &=
  \R \setminus (\M\cup\M^*\cup \m),
\end{align}
where the parameters $P=P(X)>1$ and $R=R(X)>1/\eta$ are chosen later 
(see \eqref{P-choice} and \eqref{R-choice}) as well as $\eta=\eta(X)$, that, 
as we explained before, we would like to be a small negative power of
$\max \{ p_1, p_2, p_3^k \}$ (and so of $X$). 
We have to distinguish two cases in the previous decomposition of the
real line in order to avoid a gap between the end of the major arc and the
beginning of the minor arc, where we can prove Lemma~\ref{lemma_mu}
in the form that we need: see the comments at the beginning of section
\ref{sect:inter-arc} and just before the statement of Lemma~\ref{lemma_mu}.
As we will see later in section \ref{sect:inter-arc},
we need to introduce intermediate arc only for $k \ge 5 / 2$.

The constraints on $\eta$ are in \eqref{eta-bound2}, \eqref{eta-bound}
and \eqref{eta-bound3}, according to the value of $k$.
In any case, we have
$
\I(\eta,\omega,\R)
=
\I(\eta,\omega,\M)
+
\I(\eta,\omega,\M^*)
+
\I(\eta,\omega,\m)
+
\I(\eta,\omega,\gotht)$.
We expect that $\M$ provides the main term with the right order of magnitude
without any special hypothesis on the coefficients $\lambda_j$.
It is necessary to prove that
 $\I(\eta,\omega,\M^*)$, $\I(\eta,\omega,\m)$ and 
 $\I(\eta,\omega,t)$ are $\odim{\I(\eta,\omega,\M)}$ as
$X \to +\infty$ on the particular sequence chosen:
we show that the contribution from trivial arc is ``tiny'' with
respect to the main term. The main difficulty is to estimate the minor
arc contribution; in this case we will need the full force of the
hypothesis on the coefficients $\lambda_j$ and the theory of continued
fractions.

\textbf{Remark}: from now on, anytime we use the symbol $\ll$ or $\gg$ we drop the 
dependence of the approximation from the constants $\lambda_j, \delta$ and $k$.
We use the notation $f = \infty(g)$ for $g = \odi{f}$.

\section{Lemmas}

In their original paper \cite{davenport1946indefinite} Davenport and
Heilbronn approximate directly the difference $\vert
S_k(\alpha)-T_k(\alpha)\vert $ estimating it with $\Odi{1}$.
The $L^2$-norm estimation approach (see Br\"udern, Cook \& Perelli
\cite{brudern-cook-perelli} and \cite{Languasco-Zaccagnini2016})
improves on this taking the $L^2$-norm of
$\vert S_k(\alpha)-T_k(\alpha)\vert $: this leads to the possibility
of having a wider major arc compared to the original approach.
We introduce the generalized version of the Selberg integral
\[
\mathcal{J}_k(X,h)
=
\int_{X}^{2X}
\bigl(\theta((x+h)^{1/k})-\theta(x^{1/k})-((x+h)^{1/k}-x^{1/k})\bigr)^2\, \dx x,
\]
where $\theta(x) = \sum_{p \le x} \log p$ is the usual Chebyshev function.
We have the following lemmas.

\begin{Lemma}[\cite{languasco-zaccagnini-ternary}, Theorem 3.1]
\label{norm_S-U}
Let $k\ge1$ be a real number. For $0<Y \le1/2$ we have
\[
\int_{-Y}^{Y}\vert S_k(\alpha)-U_k(\alpha)\vert ^2\, \dx \alpha
\ll
\frac{{X^{2/k-2}}\log^2 X}{Y}+Y^2X+Y^2\mathcal{J}_k\Bigl(X,\frac{1}{2Y}\Bigr).
\]
\end{Lemma}

\begin{Lemma}[\cite{languasco-zaccagnini-ternary}, Theorem 3.2]
\label{selberg}
Let $k\ge1$ be a real number and $\eps$ be an arbitrarily small positive constant. 
There exists a positive constant $c_1(\eps)$, which does not depend on $k$, such that
\[
\mathcal{J}_k(X,h)
\ll
h^2X^{2/k-1}\exp\Bigl(-c_1\Bigl(\frac{\log X}{\log \log X}\Bigr)^{1/3}\Bigr)
\]
uniformly for $X^{1-5/(6k)+\eps}\le h\le X$.
\end{Lemma}

In order to prove our crucial Lemma~\ref{Sk_quarta_tau} on the
$L^4$-norm of $S_k(\alpha)$, we need the following technical result.

\begin{Lemma}
\label{rs}
Let $\eps>0$, $k>1$ and $\gamma>0$. Let further  $\B(X^{1/k};k;\gamma)$ denote 
the number of solutions of the inequalities
\[
  \vert n_1^k+n_2^k-n_3^k-n_4^k\vert
  <
  \gamma,\qquad X^{1/k} < n_1,n_2,n_3,n_4\le 2X^{1/k}.
\]
Then
\[
\B(X^{1/k};k;\gamma)
\ll
\bigl(X^{2/k}+ \gamma X^{4/k-1} \bigr)X^{\eps}.
\]
\end{Lemma}

\begin{pf}
This is an immediate consequence of Theorem 2 of Robert \& Sargos
\cite{robert-sargos};
we just need to choose $M=X^{1/k}$, $\alpha=k$ and $\gamma=\delta M^k$ there.
\qed
\end{pf}

\begin{Lemma}
\label{Sk_quarta_tau}
Let $\eps>0$,  $\delta>0$, $k>1$, $n\in \N$ and $\tau>0$. Then we have
\[
\int_{-\tau}^{\tau}\vert S_k(\alpha)\vert ^4\, \dx \alpha
\ll
\bigl(\tau X^{2/k}+X^{4/k-1}\bigr)X^{\eps}
\quad 
\textrm{and}
\quad
\int_{n}^{n+1}\vert S_k(\alpha)\vert ^4\, \dx \alpha
\ll
\bigl( X^{2/k}+X^{4/k-1}\bigr)X^{\eps}.
\]
\end{Lemma}

\begin{pf}
A direct computation gives
\begin{align}
\notag
\int_{-\tau}^{\tau}\vert S_k(\alpha)\vert ^4\, \dx \alpha 
&
=\sum_{\delta X<p_1^k,p_2^k,p_3^k,p_4^k\le X}(\log p_1)\cdots(\log p_4)
\int_{-\tau}^{\tau}e((p_1^k+p_2^k-p_3^k-p_4^k)\alpha)\, \dx \alpha
 \\
 \notag
&\ll
(\log X)^4\sum_{\delta X<p_1^k,p_2^k,p_3^k,p_4^k\le X}
 \min\Biggl\{\tau, \frac{1}{\vert p_1^k+p_2^k-p_3^k-p_4^k\vert }\Biggr\}
\\
\notag
&\ll
(\log X)^4\sum_{\delta X<n_1^k,n_2^k,n_3^k,n_4^k\le X}
\min\Biggl\{\tau, \frac{1}{\vert n_1^k+n_2^k-n_3^k-n_4^k\vert }\Biggr\}
\\
\label{first-split-lemma-media-quarta}
&\ll 
U\tau(\log X)^4+V(\log X)^4,
\end{align}
where
\[
U:=\sum_{\substack{\delta X<n_1^k,n_2^k,n_3^k,n_4^k\le X \\ \vert n_1^k+n_2^k-n_3^k-n_4^k\vert \le 1/\tau}} 1,
\quad 
\textrm{and}
\quad
V:=\sum_{\substack{\delta X<n_1^k,n_2^k,n_3^k,n_4^k\le X \\ \vert n_1^k+n_2^k-n_3^k-n_4^k\vert > 1/\tau}}
\frac{1}{\vert n_1^k+n_2^k-n_3^k-n_4^k\vert },
\]
say.
Using Lemma \ref{rs} on $U$ we get
\begin{equation}
\label{U-estim}
U
\ll 
\B(X^{1/k};k;1/\tau)
\ll 
\Bigl(X^{2/k}+\frac1{\tau} X^{4/k-1}\Bigr)X^{\eps}.
\end{equation}
%
Concerning $V$, by a dyadic argument we get
\begin{align}
\notag
V
&\ll
\log X\Bigl(\max_{1/\tau<W\ll X}
\sum_{\substack{\delta X<n_1^k,n_2^k,n_3^k,n_4^k\le X \\ W< \vert n_1^k+n_2^k-n_3^k-n_4^k \vert \le 2W}}
  \frac1{\vert n_1^k+n_2^k-n_3^k-n_4^k\vert} \Bigr) \\
&\ll 
\log X\max_{1/\tau<W\ll X} \Bigl(\frac1W\B (X^{1/k};k;2W )\Bigr)
\notag
\ll 
 \max_{1/\tau<W\ll X}\Bigl(X^{4/k-1}+\frac{X^{2/k}}{W}\Bigr)X^{\eps}
\\
\label{V-estim}
&\ll 
(\tau X^{2/k}+X^{4/k-1})X^{\eps}.
\end{align}
Combining \eqref{first-split-lemma-media-quarta}-\eqref{V-estim}, the first part
of the lemma follows. The second part can be proved in a similar way.
\qed
\end{pf}

We need the following result in the proof of
Lemma~\ref{Lemma:bd-minor} and also when dealing with $\M^*$; see
section~\ref{sect:inter-arc}.

\begin{Lemma}
\label{tolev-pulito}
Let $\delta>0$, $k>1$, $n\in \N$ and $\tau>0$. Then
\[
\int_{-\tau}^{\tau}\vert S_k(\alpha)\vert ^2\, \dx \alpha
\ll
\bigl(\tau X^{1/k}+X^{2/k-1}\bigr) (\log X)^3
\quad
\textrm{and}
\quad
\int_{n}^{n+1}\!\!\! \vert S_k(\alpha)\vert ^2\, \dx \alpha
\ll
 X^{1/k} (\log X)^3.
\]
\end{Lemma}
\begin{pf} 
It follows directly from the proof of Lemma 7 of Tolev \cite{tolev} by letting
$c=k$ and using $X^{1/k}$  instead of $X$ there.  
We explicitly remark that the condition
$c\in (1, 15/14 )$ in Tolev's original version of this lemma
depends on other parts of his paper; in 
fact the proof of Lemma 7 of  \cite{tolev} holds for every $c>1$.
\qed
\end{pf}

We now state some other lemmas which will be mainly useful
on the minor and trivial arcs.

\begin{Lemma}[Vaughan \cite{vaughan1997hardy}, Theorem 3.1]\label{vaughan}
Let $\alpha$ be a real number and $a,q$ be positive integers satisfying $(a,q)=1$ and 
$ \vert \alpha-a/q\vert < 1/q^2$. Then
\[
S_1(\alpha)\ll\Bigl(\frac{X}{\sqrt{q}}+\sqrt{X q}+X^{4/5}\Bigr)(\log X)^4.
\]
\end{Lemma}

\begin{Lemma}\label{for-inter-arc}
Let $X^{-1} \ll \vert \alpha \vert \ll X^{-3/5}$. Then
$S_1(\alpha) \ll X^{1/2} \vert \alpha \vert^{-1/2}(\log X)^4$.
\end{Lemma}

\begin{pf}
It follows immediately from Lemma \ref{vaughan}
by choosing $q=\lfloor 1/\alpha \rfloor$ and $a=1$.
\qed
\end{pf}

\begin{Lemma}\label{harman_s1}
Let $\lambda\in \R\setminus\{0\}$, 
$X\ge Z\ge X^{4/5} (\log X)^5$ and $\vert S_1(\lambda \alpha)\vert >Z$.
Then there are coprime integers $(a,q)=1$ satisfying
\[
  1\le q
  \ll
  \Bigl(\frac{X (\log X)^4}{Z}\Bigr)^2,\qquad 
  \vert  q\lambda \alpha-a \vert 
  \ll
  \frac{X(\log X)^{10}}{Z^2}.
\]
\end{Lemma}

\begin{pf}
Let $Q$ be a parameter that we will choose later. By Dirichlet's theorem there exist 
coprime integers $(a,q)=1$ such that $1\le q \le Q$ and
$\vert q\lambda\alpha-a\vert \ll Q^{-1}\le q^{-1}$. The choice
\[
  Q = \frac{Z^2}{X(\log X)^{10}}
\]
allows us to prove the second part of the statement and to neglect
some terms in the estimations of $\vert S_1(\lambda \alpha)\vert $.
Using Lemma \ref{vaughan}, knowing that $Z\ge X^{4/5} (\log X)^5$ and
$\vert S_1(\lambda \alpha)\vert > Z$, we can rewrite the bound for
$\vert S_1(\lambda \alpha)\vert $ neglecting the term $X^{4/5}$:
\[
Z
< 
\vert S_1(\lambda \alpha)\vert 
\ll 
  (X q^{-1/2} + X^{1/2} q^{1/2}) (\log X)^4.
\]
The condition $q\le Q$ allows us to neglect the term
$X^{1/2}q^{1/2}$ and deal with small values of $q$; in fact, if
$q>X^{1/2}$ then we would have a contradiction
\[
Z
< 
\vert S_1(\lambda \alpha)\vert 
\ll 
X^{1/2}q^{1/2}(\log X)^4
\le 
X^{1/2}\frac{Z}{X^{1/2}(\log X)^5}(\log X)^4
=
\odi{Z}.
\]
Then $q \le \min \{ X^{1/2}, Q \} = X^{1/2}$, since
$Z=X^{4/5} (\log X)^5 >X^{3/4}(\log X)^5$.
Moreover, we can rewrite the inequality on $\vert S_1(\lambda \alpha)\vert $
as
\[
Z
< 
\vert S_1(\lambda \alpha)\vert 
\ll 
 X q^{-1/2}(\log X)^4
\]
and finally we get $q^{1/2} Z \ll X(\log X)^4$, which completes the proof.
\qed
\end{pf}

The optimizations in section~\ref{sect:minor-arc} depend either on
$L^2$ or on $L^4$ averages of $S_k$, according to the value of $k$;
these are provided by the following Lemmas.
For brevity, we skip the proof of the first one, remarking that it
requires Lemma~\ref{tolev-pulito}.

\begin{Lemma}[Lemma 5 of \cite{Languasco-Zaccagnini2016}]
\label{Lemma:bd-minor}
Let  $\lambda\in \R\setminus\{0\}$,
$k>1$, $0<\eta<1$, $R>1/\eta$ and $1<P<X$. We have
\[
  \int_{P/X}^{R}
    \vert S_1(\lambda \alpha) \vert^2 K_{\eta}(\alpha) \, \dx  \alpha
\ll 
  \eta X \log X
  \quad\text{and}\quad
  \int_{P/X}^{R}
    \vert S_k(\lambda  \alpha) \vert^2 K_{\eta}(\alpha) \, \dx  \alpha
\ll 
  \eta X^{1/k} (\log X)^3.
\]
\end{Lemma}

\begin{Lemma}
\label{Sk_quarta}
Let $\lambda\in \R\setminus\{0\}$, $\eps>0$, 
$k>1$, $0<\eta<1$ , $R>1/\eta$ and $1<P<X$.  Then
\[
  \int_{P/X}^{R}
\vert S_k(\lambda \alpha)\vert ^4K_{\eta}(\alpha)\, \dx \alpha
\ll 
\eta  \max\{ X^{2/k}, X^{4/k-1}\} X^{\eps}.
\]
\end{Lemma}

\begin{pf}
Using \eqref{k_eta} we obtain
\begin{equation}
\label{AB-def-lemma}
  \int_{P/X}^{R}
  \vert S_k(\lambda \alpha)\vert ^4K_{\eta}(\alpha)\, \dx \alpha
\ll
\eta^2\int_{P/X}^{1/\eta}\vert S_k(\lambda \alpha)\vert ^4\, \dx \alpha
+
\int_{1/\eta}^{R}\vert S_k(\lambda \alpha)\vert ^4\frac{\dx \alpha}{\alpha^2}
=
A+ B,
\end{equation}
say.
By  Lemma \ref{Sk_quarta_tau}, we immediately get
\begin{equation}
\label{A-estim-lemma}
A
\ll 
\eta^2 \int_{- \vert \lambda\vert /\eta}^{ \vert \lambda\vert /\eta} \vert S_k(\alpha)\vert ^4\, \dx \alpha
\ll 
\eta \max\{ X^{2/k}, \eta X^{4/k-1} \} X^{\eps}.
\end{equation}
Moreover, again by  Lemma \ref{Sk_quarta_tau},  we have that
\begin{align}
\notag
B
&\ll
\int_{\vert \lambda\vert /\eta}^{+\infty}\vert S_k(\alpha)\vert ^4\frac{\dx \alpha}{\alpha^2}
\ll  
\sum_{n\ge \vert \lambda\vert /\eta}\frac{1}{(n-1)^2}\int_{n-1}^n\vert S_k(\alpha)\vert ^4\, \dx \alpha
\\
\label{B-estim-lemma}
&
\ll 
\eta \max\{ X^{2/k}, X^{4/k-1} \} X^{\eps}.
\end{align}
Combining \eqref{AB-def-lemma}-\eqref{B-estim-lemma} and using $0<\eta<1$, the lemma follows.
\qed
\end{pf}

As we remarked in the introduction, stronger bounds are now available
for larger integral $k$, but they are not useful for our purpose.
The next Lemma provides the additional information that enables us to
give a non-trivial result also when $k = 3$.

\begin{Lemma}
\label{hua-app}
Let $\lambda\in \R\setminus\{0\}$, $\eps>0$,
$0<\eta<1$, $R>1/\eta$ and $1<P<X$. Then 
\[
  \int_{P/X}^{R}
\vert S_3(\lambda \alpha)\vert ^8 K_{\eta}(\alpha)\, \dx \alpha
\ll 
\eta    X^{5/3+\eps}.
\]
\end{Lemma}

\begin{pf}
Inserting  Hua's estimate in \cite{Hua1938}, i.e.
$
\int_0^1 \vert S_3(\alpha)\vert ^8 \ \dx  \alpha\ll X^{5/3+\eps},
$
in the body of the proof of Lemma \ref{Sk_quarta}
and exploiting the periodicity of $S_3(\alpha)$,
the result follows immediately.
\qed
\end{pf}
 
Another lemma on the minor arc is inserted in the body of section
\ref{sect:minor-arc}.

\section{The major arc}

We recall the definitions in \eqref{def-arcs-k-basso} and
\eqref{def-arcs-k-alto}.
The major arc computation is the same as in \cite{Languasco-Zaccagnini2016}:
\begin{align*}
\I(\eta,\omega,\M)
&=
\int_{\M}S_1(\lambda_1\alpha)S_1(\lambda_2\alpha)S_k(\lambda_3\alpha)
K_{\eta}(\alpha)e(-\omega\alpha)\, \dx \alpha
\\
&=
\int_{\M}T_1(\lambda_1\alpha)T_1(\lambda_2\alpha)T_k(\lambda_3\alpha)
K_{\eta}(\alpha)e(-\omega\alpha)\, \dx \alpha
\\
&\hskip1cm
+\int_{\M}(S_1(\lambda_1\alpha)-T_1(\lambda_1\alpha))T_1(\lambda_2\alpha)T_k(\lambda_3\alpha)
K_{\eta}(\alpha)e(-\omega\alpha)\, \dx \alpha
\\
&\hskip1cm
+\int_{\M}S_1(\lambda_1\alpha)(S_1(\lambda_2\alpha)-T_1(\lambda_2\alpha))T_k(\lambda_3\alpha)K_{\eta}(\alpha)e(-\omega\alpha)\, \dx \alpha
\\
&
\hskip1cm
+\int_{\M}S_1(\lambda_1\alpha)S_1(\lambda_2\alpha)(S_k(\lambda_3\alpha)-T_k(\lambda_3\alpha))K_{\eta}(\alpha)e(-\omega\alpha)\, \dx \alpha
\\
&=J_1+J_2+J_3+J_4,
\end{align*}
say.

\subsection{Main term: lower bound for \texorpdfstring{$J_1$}{J1}}

As the reader might expect the main term is given by the summand $J_1$.

Let $H(\alpha)=T_1(\lambda_1\alpha)T_1(\lambda_2\alpha) T_k(\lambda_3\alpha)K_{\eta}(\alpha)e(-\omega\alpha)$ so that
\[
J_1
=
\int_{\R}H(\alpha)\, \dx \alpha
+
\Odig{\int_{P/X}^{+\infty}\vert H(\alpha)\vert \, \dx \alpha}.
\]
Using \eqref{k_eta} and \eqref{stima_tk}, we get
\[
\int_{P/X}^{+\infty}\vert H(\alpha)\vert \, \dx \alpha
\ll  \eta^2X^{1/k-1}
\int_{P/X}^{+\infty}\frac{\dx \alpha}{\alpha^3}
\ll  
\eta^2\frac{X^{1+1/k}}{P^2}
=
\odim{\eta^2X^{1+1/k}},
\]
provided that $P\rightarrow+\infty$.
Let now $D=[\delta X,X]^2\times[(\delta X)^{1/k},X^{1/k}]$. We obtain
\begin{align*}
\int_{\R}H(\alpha)\, \dx \alpha
&=
\iiint_D
\int_{\R}
e((\lambda_1t_1+\lambda_2t_{2}+\lambda_3t_3^k-\omega)\alpha)
K_{\eta}(\alpha)\, \dx \alpha\,
\dx  t_1\dx  t_2\dx  t_3
\\
&=
\iiint_D
\max\{ 0, \eta-\vert \lambda_1t_1+\lambda_2t_2+\lambda_3t_3^k-\omega)\vert \} \, 
\dx  t_1\dx  t_2\dx  t_3.
\end{align*}

Apart from trivial changes of sign, there are essentially two cases:
\begin{enumerate}
\item $\lambda_1>0$, $\lambda_2>0$, $\lambda_3<0$ 
\item $\lambda_1>0$, $\lambda_2<0$, $\lambda_3<0$.
\end{enumerate}

We deal with the first one.
We warn the reader that here it may be necessary to adjust the value
of $\delta$ in order to guarantee the necessary set inclusions.
After a suitable change of variables, letting
$D'=[\delta X,(1-\delta)X]^3$, we find that
\begin{align*}
  \int_{\R}H(\alpha)\, \dx \alpha
  &\gg
  \iiint_{D'}
  \max\{ 0, \eta-\vert \lambda_1u_1+\lambda_2u_2+\lambda_3u_3)\vert \} \, 
    u_3^{1/k-1}\ 
  \dx  u_1\dx  u_2\dx  u_3 \\
  &\gg
  X^{1/k-1}
  \iiint_{D^{\prime}}
  \max\{ 0, \eta-\vert \lambda_1u_1+\lambda_2u_2+\lambda_3u_3)\vert \} \, 
   \dx  u_1\dx  u_2\dx  u_3.
\end{align*}
Apart from sign, the computation is essentially symmetrical with
respect to the coefficients $\lambda_j$: we assume, as we may, that
$| \lambda_3| \ge \max\{ \lambda_1, \lambda_2 \}$, the other cases
being similar.
Now, for $j=1,2$ let $a_j=\dfrac{2\delta\vert \lambda_3\vert }{\vert \lambda_j\vert }$, $b_j=\dfrac32 a_j$ and 
$\mathcal{D}_j=[a_jX,b_jX]$; if $u_j\in\mathcal{D}_j$ for $j=1,2$ then
\[
\lambda_1u_1+\lambda_2u_2\in\left[4\vert \lambda_3\vert \delta X,6\vert \lambda_3\vert \delta X\right]
\]
so that, for every choice of $(u_1,u_2)$ the interval $[a,b]$ with endpoints 
$\pm\eta/\vert \lambda_3\vert +(\lambda_1 u_1+\lambda_2 u_2)/\vert \lambda_3\vert $ is contained 
in $[\delta X,(1-\delta)X]$. In other words, for $u_3\in[a,b]$ the values of
 $\lambda_1u_1+\lambda_2u_2+\lambda_3u_3$ cover the whole interval $[-\eta,\eta]$.
Hence for any $(u_1,u_2)\in\mathcal{D}_1\times\mathcal{D}_2$ we have
\[
\int_{\delta X}^{(1-\delta)X}
\max\{ 0, \eta-\vert \lambda_1u_1+\lambda_2u_2+\lambda_3u_3\vert \} \, \dx u_3
=
\vert \lambda_3\vert ^{-1}
\int_{-\eta}^{\eta}\max\{ 0,\eta-\vert u\vert \} \, \dx u\gg\eta^2.
\]
Summing up, we get
\[
J_1
\gg
\eta^2X^{1/k-1}
\iint_{\mathcal{D}_1\times\mathcal{D}_2}\ \dx  u_1\dx  u_2
\gg
\eta^2X^{1/k-1}X^2
=
\eta^2X^{1+1/k},
\]
which is the expected lower bound.

\subsection{Bound for \texorpdfstring{$J_2$, $J_3$ and $J_4$}{J2, J3 and J4}}

The computations for $J_2$ and $J_3$ are similar to and simpler than the
corresponding one for $J_4$; moreover the most restrictive condition on
$P$ arises from $J_4$; hence we will skip the computation for both
$J_2$ and $J_3$.
Using the triangle inequality and \eqref{k_eta},
\begin{align*}
J_4
&\ll 
\eta^2\int_{\M}\vert S_1(\lambda_1\alpha)\vert 
\vert S_1(\lambda_2\alpha)\vert 
\vert  S_k(\lambda_3\alpha)-T_k(\lambda_3\alpha) \vert \, 
\dx \alpha
\\ 
&\le 
\eta^2\int_{\M}\vert S_1(\lambda_1\alpha)\vert 
\vert S_1(\lambda_2\alpha)\vert  
\vert S_k(\lambda_3\alpha)-U_k(\lambda_3\alpha)\vert \, 
\dx \alpha
\\
&\hskip 1cm +
\eta^2\int_{\M}\vert S_1(\lambda_1\alpha)\vert 
\vert S_1(\lambda_2\alpha)\vert 
\vert U_k(\lambda_3\alpha)-T_k(\lambda_3\alpha)\vert \, 
\dx \alpha
\\
&=
\eta^2(A_4+B_4),
\end{align*}
say, where $U_k(\lambda_3\alpha)$ and $T_k(\lambda_3\alpha)$ are defined in  \eqref{U_k}.
Using the Cauchy-Schwarz inequality, Lemmas
\ref{norm_S-U}-\ref{selberg} and trivial bounds yields, for any fixed
$A>0$,
\begin{align*}
A_4
& \ll  
X
\Bigl(
\int_{\M}\vert S_1(\lambda_1\alpha)\vert ^2\, \dx \alpha
\Bigr)^{1/2}
\Bigl(
\int_{\M}\vert S_k(\lambda_3\alpha)-U_k(\lambda_3\alpha)\vert ^2\, \dx \alpha
\Bigr)^{1/2}
\\
& \ll  
X^{1+1/k}(\log X)^{(1-A)/2}
=
\odim{X^{1+1/k}}
\end{align*}
as long as $A>1$, provided that $P\le X^{5/(6k)-\eps}$.
Using again the Cauchy-Schwarz inequality, \eqref{t-u} and trivial bounds,
we see that
\begin{align*}
B_4
& \ll 
\int_0^{1/X}
\vert S_1(\lambda_1\alpha)\vert 
\vert S_1(\lambda_2\alpha)\vert \, 
\dx \alpha
+
X\int_{1/X}^{P/X}
\alpha
\vert S_1(\lambda_1\alpha)\vert 
\vert S_1(\lambda_2\alpha)\vert \, 
\dx \alpha
\\
& \ll  
X+
P\Bigl(
\int_{1/X}^{P/X}\vert S_1(\lambda_1\alpha)\vert ^2\, \dx \alpha 
\int_{1/X}^{P/X}\vert S_1(\lambda_2\alpha)\vert ^2\, \dx \alpha
\Bigr)^{1/2}
\ll 
PX\log X.
\end{align*}

Taking $P=\odim{X^{1/k}(\log X)^{-1}}$ we get $\eta^2 B_4=\odim{\eta^2 X^{1+1/k}}$. 
We may therefore choose
\begin{equation}
\label{P-choice}
P=X^{5/(6k)-\eps}.
\end{equation}

\section{The trivial arc}

We recall that the trivial arc is defined in \eqref{def-arcs-k-basso}
and \eqref{def-arcs-k-alto}.
Using the Cauchy-Schwarz inequality and \eqref{stima_tk},
we see that
\begin{align*}
\vert \I(\eta,\omega,\gotht)\vert 
&
\ll \int_{R}^{+\infty}
\vert S_1(\lambda_1\alpha)S_1(\lambda_2\alpha)S_k(\lambda_3\alpha)\vert
K_{\eta}(\alpha) \, 
\dx \alpha
\\
&\ll 
X^{1/k}
\Bigl(\int_{R}^{+\infty}\vert S_1(\lambda_1\alpha)\vert ^2K_{\eta}(\alpha)\, \dx \alpha\Bigr)^{1/2}
\Bigl(\int_{R}^{+\infty}\vert S_1(\lambda_2\alpha)\vert ^2K_{\eta}(\alpha)\, \dx \alpha\Bigr)^{1/2}
\\
& \ll 
X^{1/k}C_1^{1/2}C_2^{1/2},
\end{align*}
say.
Using the PNT and the periodicity of $S_1(\alpha)$, for every $j=1,2$
we have that
\[
C_j
=
\int_{R}^{+\infty}\!\!\!\! \vert S_1(\lambda_j\alpha)\vert ^2 \frac{\dx \alpha}{\alpha^2}
\ll
\int_{\vert \lambda_j\vert R}^{+\infty} \vert S_1(\alpha)\vert ^2  \frac{\dx \alpha}{\alpha^2}
\ll 
\sum_{n\ge\vert \lambda_j\vert  R}\frac{1}{(n-1)^2}\int_{n-1}^n \!\!  \vert S_1(\alpha)\vert ^2\, \dx \alpha
\ll
\frac{X\log X}{\vert \lambda_j\vert  R}.
\]
Hence, recalling that $\vert \I(\eta,\omega, \gotht)\vert $ has to be $o (\eta^2X^{1+1/k} )$,
  the choice
\begin{equation}
\label{R-choice}
R=\eta^{-2} (\log X)^{3/2}
\end{equation}
is admissible.

\section{The intermediate arc: \texorpdfstring{$5/2 \le k \le3$}{large k}}
\label{sect:inter-arc}

In section \ref{sect:minor-arc} we apply Harman's technique to the
minor arc, using Lemma~\ref{harman_s1} as the starting point.
We remark that in the course of the proof of Lemma~\ref{lemma_mu} it
is crucial that both the integers $a_1$ and $a_2$ appearing in
\eqref{misura3} below do not vanish; in fact, if $a_1 = 0$, say, then
$\alpha$ is very small ($\alpha \ll X^{-2/3}$) and, according to our
definitions above, it belongs to $\M \cup \M^*$.

For small $k$ we do not need an intermediate arc, because the major
arc is wide enough to rule out the possibility that $a_1 a_2 = 0$ for
$\alpha \in \m$.
For larger values of $k$, the constraint \eqref{P-choice} implies that
there is a gap between the major arc and the minor arc which we need
to fill: see the definition in \eqref{def-arcs-k-alto}.
Using the intermediate arc $\M^*$, we are able to cover more than
needed.

Let $5/2 \le k\le 3$: we now show that the contribution of $\M^*$ is
negligible.
Using \eqref{k_eta}, Lemma \ref{for-inter-arc}, the Cauchy-Schwarz
inequality and \eqref{P-choice} we get
\begin{align*}
\I(\eta,\omega,\M^*)&
\ll 
\eta^2 
\int_{P/X}^{X^{-3/5}}
\vert S_1(\lambda_1\alpha)\vert 
\vert S_1(\lambda_2\alpha)\vert 
\vert S_k(\lambda_3\alpha)\vert 
\, \dx \alpha \\
&\ll
\eta^2  X(\log X)^8 
\int_{P/X}^{X^{-3/5}}
\vert S_k(\lambda_3\alpha)\vert \, \frac{\dx \alpha}{\alpha}
\\ & 
\ll  
\eta^2 X (\log X)^8
\Bigl(\int_{-X^{-3/5}}^{X^{-3/5}}\vert S_k(\lambda_3\alpha)\vert ^2\, \dx \alpha\Bigr)^{1/2}
\Bigl(\int_{P/X}^{X^{-3/5}} \frac{\dx \alpha}{\alpha^2}\Bigr)^{1/2}
\\
& 
\ll 
\eta^2 X (X^{1/k-3/5} )^{1/2} (X^{1-5/(6k)} )^{1/2} X^{\eps}
\ll 
\eta^2 X^{6/5+1/(12k)+\eps},
\end{align*}
where we also used Lemma \ref{tolev-pulito} with $\tau=X^{-3/5}$ and the
fact that $k \ge 5 / 2$.
The last estimate is $\odim{\eta^2X^{1+1/k}}$ for every
$5/2 \le k<55/12$.

\section{The minor arc}
\label{sect:minor-arc}

Here we use Harman's technique as described in \cite{harman1991general}.
The minor arc $\m$ is defined in \eqref{def-arcs-k-basso} and
\eqref{def-arcs-k-alto}, according to the value of $k$.
In view of using Lemma \ref{harman_s1}, we now split $\m$ into subsets
$\m_1$, $\m_2$ and $\m^* = \m \setminus (\m_1\cup \m_2)$, where
\[
  \m_i
  =
  \{\alpha\in \m \colon \vert S_1(\lambda_i\alpha)\vert \le X^{5/6} (\log X)^5\}
  \quad \text{for}\ i=1,2.
\]
In order to obtain the optimization, we chose to split the range for
$k$ into two intervals in which to take advantage of the $L^2$-norm of
$S_k(\alpha)$ in one case (Lemma \ref{Lemma:bd-minor}) and the
$L^4$-norm of $S_k(\alpha)$ in the other one (Lemma \ref{Sk_quarta}).
The same choice will be made later in the discussion of the arc $\m^*$.
We will see that it is not possible to split the minor arc in another
way in order to get a better result, in the present state of knowledge
on exponential sums.

\subsection{Bounds on \texorpdfstring{$\m_1\cup \m_2$}{m1 and m2}}
Using  H\"older's inequality and Lemma \ref{Lemma:bd-minor}, 
for $1<k\le 6/5$ we obtain
\begin{align}
\nonumber
\vert \I(\eta,\omega,\m_i)\vert 
& \ll
\int_{\m_i}\vert S_1(\lambda_1\alpha)\vert 
\vert S_1(\lambda_2\alpha)\vert 
\vert S_k(\lambda_3\alpha)\vert  K_{\eta}(\alpha) \, \dx \alpha
\\ 
\nonumber
&\ll 
\Bigl(\max_{\alpha\in \m_i}{\vert S_1(\lambda_1\alpha)\vert}\Bigr)
\Bigl(\int_{\m_i}\vert S_1(\lambda_2\alpha)\vert ^2 K_{\eta}(\alpha) \, \dx \alpha)\Bigr)^{1/2}
\\
\nonumber
&
\hskip1cm
\times
\Bigl(\int_{\m_i}\vert S_k(\lambda_3\alpha)\vert ^2 K_{\eta}(\alpha) \, \dx \alpha)\Bigr)^{1/2}
\\
\nonumber
&\ll
X^{5/6} (\log X)^5 (\eta X\log X)^{1/2}
(\eta X^{1/k}(\log X)^3)^{1/2}
\\
\label{cond_m2}
&
\ll
\eta X^{4/3+1/(2k)+\eps}.
\end{align}
The estimate in \eqref{cond_m2} should be $o (\eta^2X^{1+1/k} )$;
hence this leads to the constraint
\begin{equation}
\label{eta-bound2}
\eta=\infty (X^{1/3-1/(2k)+\eps} ),
\end{equation}
where $f = \infty(g)$ means $g = \odi{f}$.

Using H\"older's inequality and Lemmas \ref{Lemma:bd-minor} and
\ref{Sk_quarta}, for $6/5<k<3$ we obtain
\begin{align}
\nonumber
\vert \I(\eta,\omega,\m_i)\vert 
& \ll
\int_{\m_i}\vert S_1(\lambda_1\alpha)\vert 
\vert S_1(\lambda_2\alpha)\vert 
\vert S_k(\lambda_3\alpha)\vert  K_{\eta}(\alpha) \, \dx \alpha
\\ 
\nonumber
&\ll 
\Bigl(\max_{\alpha\in \m_i}{\vert S_1(\lambda_1\alpha)\vert ^{1/2}}\Bigr)
\Bigl(\int_{\m_i}\vert S_1(\lambda_1\alpha)\vert ^2 K_{\eta}(\alpha) \, \dx \alpha)\Bigr)^{1/4}
\\
\nonumber
&
\hskip1cm
\times
\Bigl(\int_{\m_i}\vert S_k(\lambda_3\alpha)\vert ^4 K_{\eta}(\alpha) \, \dx \alpha)\Bigr)^{1/4}
\Bigl(\int_{\m_i}\vert S_1(\lambda_2\alpha)\vert ^2 K_{\eta}(\alpha) \, \dx \alpha)\Bigr)^{1/2}
\\
\nonumber
&\ll
X^{5/12} (\log X)^{5 / 2} (\eta X\log X)^{1/4}
(\eta \max\{ X^{2/k}, X^{4/k-1} \})^{1/4}(\eta X\log X)^{1/2}
\\
\label{cond_m1}
&
\ll
\eta \max \bigl\{ X^{7/6+1/(2k)}, X^{11/12+1/k} \bigr\} X^{\eps}.
\end{align}
The estimate in \eqref{cond_m1} should be $o (\eta^2X^{1+1/k} )$;
hence this leads to
\begin{equation}
\label{eta-bound}
  \eta
  =
  \infty \bigl( \max\{ X^{1/6-1/(2k)+\eps}, X^{-1/12+\eps} \} \bigr).
\end{equation}
If $k=3$ we use Lemmas \ref{Lemma:bd-minor} and \ref{hua-app} thus getting
\begin{align*}
\vert \I(\eta,\omega, \m_i)\vert 
&\ll
\int_{\m_i}\vert S_1(\lambda_1\alpha)\vert 
\vert S_1(\lambda_2\alpha)\vert 
\vert S_3(\lambda_3\alpha)\vert  K_{\eta}(\alpha) \, \dx \alpha
\\ 
&\ll 
\Bigl(\max_{\alpha\in \m_i}{\vert S_1(\lambda_1\alpha)\vert ^{1/4}}\Bigr)
\Bigl(\int_{\m_i}\vert S_1(\lambda_1\alpha)\vert ^2 K_{\eta}(\alpha) \, \dx \alpha)\Bigr)^{3/8}
\\
& \hskip1cm \times 
\Bigl(\int_{\m_i}\vert S_3(\lambda_3\alpha)\vert ^8 K_{\eta}(\alpha) \, \dx \alpha)\Bigr)^{1/8}
\Bigl(\int_{\m_i}\vert S_1(\lambda_2\alpha)\vert ^2 K_{\eta}(\alpha) \, \dx \alpha)\Bigr)^{1/2}
\\
&\ll 
\eta X^{31/24+\eps}.
\end{align*}
This bound leads to the constraint 
\begin{equation}
\label{eta-bound3}
  \eta
  =
  \infty \bigl( X^{-1/24 + \eps} \bigr),
\end{equation}
which justifies the last line of \eqref{psi-value}.

\subsection{Bound on \texorpdfstring{$\m^*$}{m*}}
We recall our definitions in \eqref{def-arcs-k-basso} and
\eqref{def-arcs-k-alto}.
It remains to discuss the set $\m^*$ where the following bounds hold
simultaneously
\[
\vert S_1(\lambda_1\alpha)\vert > X^{5/6} (\log X)^5,
\quad 
\vert S_1(\lambda_2\alpha)\vert > X^{5/6} (\log X)^5,
\quad 
T \le \vert \alpha\vert \le \eta^{-2}(\log X)^{3/2} = R,
\]
where $T = P / X = X^{5/(6k) - 1 - \eps}$ by our choice in \eqref{P-choice}
if $k < 5 / 2$, and $T = X^{- 3 / 5}$ otherwise. 
Using a dyadic dissection, we split  $\m^*$ into disjoint sets
$E(Z_1,Z_2,y)$ in which, for $\alpha\in E(Z_1,Z_2,y)$, we have
\[
  Z_i<\vert S_1(\lambda_i\alpha)\vert \le 2Z_i,
  \qquad
  y<\vert \alpha\vert \le 2y,
\]
where $Z_i=2^{k_i}X^{5/6} (\log X)^5$ and $y=2^{k_3}X^{5/(6k)-1-\eps}$ for some 
non-negative integers $k_1,k_2,k_3$.

It follows that the number of disjoint sets is, at most, $\ll(\log X)^3$.
Let us write $\A$ as a shorthand for the set $E(Z_1,Z_2,y)$.
We need an upper bound for the Lebesgue measure of $\A$.
In the following Lemma, it is crucial that both the integers $a_1$ and
$a_2$ appearing in \eqref{misura3} below do not vanish; in fact, if
$a_1 = 0$, say, then $q_1 = 1$ and $\alpha$ is so small that it can
not belong to $\m$.
If $k$ is large, we treat the range $[P/X,X^{-3/5}]$ and its
symmetrical by means of the argument in section \ref{sect:inter-arc}:
this is needed because, in this case, the inequalities \eqref{misura3}
below do not rule out the possibility that $a_1a_2=0$, unless
$\vert\alpha\vert$ is large enough.

\begin{Lemma}
\label{lemma_mu}
Let $\eps > 0$.
We have that $\mu(\A)\ll yX^{8/3+\eps}Z_1^{-2}Z_2^{-2}$,
where $\mu(\cdot)$ denotes the Lebesgue measure.
\end{Lemma}

\begin{pf}
If $\alpha\in\A$, by Lemma \ref{harman_s1}
there are coprime integers $(a_1,q_1)$ and $(a_2,q_2)$ such that
\begin{equation}
\label{misura3}
 1\le q_i\ll\Bigl(\frac{X(\log X)^4}{Z_i}\Bigr)^2,
\qquad 
\vert q_i\lambda_i\alpha-a_i\vert 
\ll
 \frac{X(\log X)^{10}}{Z_i^2}.
\end{equation}

We remark that $a_1a_2\neq0$ otherwise we would have $\alpha\in\M \cup \M^*$.
In fact, if $a_1a_2=0$, recalling the definitions of $Z_i$ and \eqref{misura3},
$\alpha \ll q_i^{-1} X(\log X)^{10}Z_i^{-2} \ll X^{-2/3}$.

Now, we can further split $\m^*$ into sets $I = I(Z_1,Z_2,y,Q_1,Q_2)$ where, on each set,
 $Q_j\le q_j\le 2Q_j$. Note that $a_i$ and $q_i$ are uniquely determined by $\alpha$; in the opposite direction, for a given quadruple $a_1$, $q_1$, $a_2$, $q_2$, the inequalities \eqref{misura3} define an interval of $\alpha$ of length
\[
  \ll
  \min\Biggl\{ \frac{
X(\log X)^{10}}{Q_1Z_1^2},\frac{
X(\log X)^{10}}{Q_2Z_2^2} \Biggr\} 
  \ll
  \frac{X(\log X)^{10}}{Q_1^{1/2}Q_2^{1/2}Z_1Z_2},
\]
by taking the geometric mean.

Now we need a lower bound for $Q_1 Q_2 $: by \eqref{misura3} we obtain
\begin{align*}
\Bigl\vert a_2q_1\frac{\lambda_1}{\lambda_2}-a_1q_2\Bigr\vert 
&
=
\Bigl\vert 
\frac{a_2}{\lambda_2\alpha}(q_1\lambda_1\alpha-a_1)-\frac{a_1}{\lambda_2\alpha}(q_2\lambda_2\alpha-a_2)
\Bigr\vert  
\\
&\ll 
q_2\vert q_1\lambda_1\alpha-a_1\vert +q_1\vert q_2\lambda_2\alpha-a_2\vert  
\\
& \ll
Q_2 \frac{X(\log X)^{10}}{Z_1^2}
+
Q_1\frac{X(\log X)^{10}}{Z_2^2}.
\end{align*}

Recalling that $Q_i\ll  ( X(\log X)^{4}/Z_i )^2$ and that
$Z_i \gg X^{5/6} (\log X)^5$, we have 
\begin{equation}
\label{harman_stima3}
\Bigl\vert a_2q_1\frac{\lambda_1}{\lambda_2}-a_1q_2\Bigr\vert 
\ll
\Bigl(\frac{X(\log X)^{4}}{X^{5/6} (\log X)^5}\Bigr)^2
\Bigl(\frac{X^{1/ 2}(\log X)^5}{X^{5/6} (\log X)^5}\Bigr)^2 
\ll 
  X^{-1/3} (\log X)^{-2}
  <\frac{1}{4q}.
\end{equation}

We recall that $q=X^{1/3}$ is a denominator of a convergent of $\lambda_1/\lambda_2$. Hence by \eqref{harman_stima3}, Legendre's law of best approximation 
for continued fractions implies that $\vert a_2q_1\vert\ge q$ and by the same token, for any pair $\alpha$, $\alpha'$ having distinct associated products $a_2q_1$,
\[
\vert a_2(\alpha)q_1(\alpha)-a_2(\alpha')q_1(\alpha')\vert \ge q;
\]
thus, by the pigeon-hole principle, there is at most one value of $a_2q_1$ in the 
interval $[rq,(r+1)q)$ for any positive integer $r$. Furthermore $a_2q_1$ determines $a_2$ and 
$q_1$ to within $X^{\eps/2}$ possibilities (from the bound for the divisor function)
 and consequently also $a_2q_1$ determines $a_1$ and $q_2$ to within $X^{\eps/2}$ 
 possibilities from \eqref{harman_stima3}.

Hence we got a lower bound for $q_1q_2$, since, using $Q_j\le q_j\le 2Q_j$, we get
\[
q_1q_2=a_2q_1\frac{q_2}{a_2}
\gg
\frac{rq}{\vert \alpha\vert }
\gg 
rqy^{-1}.
\]
for the quadruple under consideration.

As a consequence we obtain that the total length of the part of $I(Z_1,Z_2,y,Q_1,Q_2)$ with $a_2q_1\in[rq,(r+1)q)$ is
\[
\ll
  X^{1 + \eps/2} (\log X)^{10}{Z_1^{-1}Z_2^{-1}}r^{-1/2}q^{-1/2}y^{1/2}.
\]

Now we need a bound for $r$: since $a_2 q_1 \in [r q, (r + 1),q)$, we have
\[
rq\le\vert a_2q_1\vert\ll q_1q_2\vert \alpha\vert 
\ll
y \Bigl(\frac{X(\log X)^{4}}{Z_1}\Bigr)^2\Bigl(\frac{X(\log X)^{4}}{Z_2}\Bigr)^2 
\ll
  \frac{y X^{4} (\log X)^{16}}{Z_1^2 Z_2^2}
\]
and hence we get
\[
r\ll q^{-1}yX^{4} (\log X)^{16}Z_1^{-2}Z_2^{-2}.
\]
Next, we sum on every interval to get an upper bound for the measure
of $\A$: we get
\[
  \mu(\A)
  \ll
  \frac{X^{1 + \eps/2} y^{1/2} (\log X)^{10}}{Z_1 Z_2 q^{1/2}}
  \sum_{1\le r\ll q^{-1}yX^{4} (\log X)^{16}Z_1^{-2}Z_2^{-2}}r^{-1/2}.
\]
Standard estimates imply that the sum on the right is
$\ll (q^{-1}yX^{4} (\log X)^{16}Z_1^{-2}Z_2^{-2})^{1/2}$, and recalling that
$q = X^{1/3}$ we can finally write
\[
  \mu(\A)
  \ll 
  y X^{3 + \eps/2} (\log X)^{18}Z_1^{-2}Z_2^{-2}q^{-1}
  \ll
  y X^{8/3 + \eps}Z_1^{-2}Z_2^{-2}.
\]
This proves the lemma.
\qed
\end{pf}

\section{Conclusion}
Here we finally justify the choice of the function $\psi$ in the statement
of the main Theorem. 
Using Lemmas \ref{Lemma:bd-minor}-\ref{Sk_quarta}-\ref{lemma_mu} we are
now able to estimate $\I(\eta,\omega,\A)$ for $1 < k \le 3$. 
For $k \ge \frac52$, we also need the result in section \ref{sect:inter-arc}.

If $1<k\le 6/5$ we proceed as follows:
\begin{align*}
\vert \I(\eta,\omega, \A) \vert
&\ll
\int_{\A}
\vert S_1(\lambda_1\alpha)\vert 
\vert S_1(\lambda_2\alpha)\vert 
\vert S_k(\lambda_3\alpha)\vert  
K_{\eta}(\alpha) \, \dx \alpha
\\
&\ll
 \Bigl(
 \int_{\A}
 \vert S_1(\lambda_1\alpha)S_1(\lambda_2\alpha)\vert ^2
  K_{\eta}(\alpha) \, \dx \alpha
  \Bigr)^{1/2}
 \Bigl(
 \int_{\A}
 \vert S_k(\lambda_3\alpha)\vert ^2 
 K_{\eta}(\alpha) \, \dx \alpha
 \Bigr)^{1/2}
 \\
&\ll 
\bigl(
\min\bigl\{\eta^2, y^{-2}\bigr\}
\bigr)^{1/2}
\bigl(
(Z_1Z_2)^2\mu(\A)
\bigr)^{1/2}
\bigl(
\eta  X^{1/k+\eps}
\bigr)^{1/2}
\\
&\ll 
\bigl(
\min\bigl\{\eta^2, y^{-2}\bigr\}
\bigr)^{1/2}
Z_1Z_2
\big(
yX^{8/3 + \eps} Z_1^{-2}Z_2^{-2}
\big)^{1/2}
\eta^{1/2}
X^{1/(2k)+\eps/2}
\\
&\ll 
\eta X^{4/3+1/(2k)+\eps}.
\end{align*}
Hence we need $\eta = \infty \bigl( X^{1/3-1/(2k)+\eps} \bigr)$, which
is  the same condition we got in \eqref{eta-bound2}.

If $6/5<k<3$, 
\begin{align*}
\vert \I(\eta,\omega, \A) \vert
&\ll
\int_{\A}
\vert S_1(\lambda_1\alpha)\vert 
\vert S_1(\lambda_2\alpha)\vert 
\vert S_k(\lambda_3\alpha)\vert  
K_{\eta}(\alpha) \, \dx \alpha
\\
&\ll
 \Bigl(
 \int_{\A}
 \vert S_1(\lambda_1\alpha)S_1(\lambda_2\alpha)\vert ^{4/3}
  K_{\eta}(\alpha) \, \dx \alpha
  \Bigr)^{3/4}
 \Bigl(
 \int_{\A}
 \vert S_k(\lambda_3\alpha)\vert ^4 
 K_{\eta}(\alpha) \, \dx \alpha
 \Bigr)^{1/4}
 \\
&\ll 
\bigl(
\min\bigl\{\eta^2, y^{-2}\bigr\}
\bigr)^{3/4}
\bigl(
(Z_1Z_2)^{4/3}\mu(\A)
\bigr)^{3/4}
\bigl(
\eta \max\{ X^{2/k}, X^{4/k-1}\} X^{\eps}
\bigr)^{1/4}
\\
&\ll 
\bigl(
\min\bigl\{\eta^2, y^{-2}\bigr\}
\bigr)^{3/4}
Z_1Z_2
\big(
yX^{8/3 + \eps} Z_1^{-2}Z_2^{-2}
\big)^{3/4}
\eta^{1/4}
\max\{ X^{1/(2k)}, X^{1/k-1/4} \} X^{\eps / 4}
\\
&\ll 
\eta Z_1^{-1/2}Z_2^{-1/2}X^{2+\eps}
\max\{ X^{1/(2k)}, X^{1/k-1/4} \}
\\
&\ll 
\eta 
\max\{ X^{7/6+1/(2k)}, X^{11/12+1/k} \} X^{\eps}.
\end{align*}
Hence we need
$\eta = \infty \bigl( \max\{ X^{1/6-1/(2k)+\eps}, X^{-1/12+\eps} \} \bigr)$,
which is the same condition we got in \eqref{eta-bound}.
If $k = 3$, using Lemmas \ref{hua-app} and \ref{lemma_mu} we obtain
\begin{align*}
\vert \I(\eta,\omega,\A) \vert
&\ll
\int_{\A}
\vert S_1(\lambda_1\alpha)\vert 
\vert S_1(\lambda_2\alpha)\vert 
\vert S_3(\lambda_3\alpha)\vert  
K_{\eta}(\alpha) \, \dx \alpha
\\
&\ll 
\Bigl(
\int_{\A}
\vert S_1(\lambda_1\alpha)S_1(\lambda_2\alpha)\vert ^{8/7} 
K_{\eta}(\alpha) \, \dx \alpha
\Bigr)^{7/8}
\Bigl(
\int_{\A}
\vert S_3(\lambda_3\alpha)\vert ^8 
K_{\eta}(\alpha) \, \dx \alpha
\Bigr)^{1/8}
\\
&\ll 
\eta Z_1^{-3/4}Z_2^{-3/4}X^{7/3+5/24+\eps}
\ll
\eta X^{31/24+\eps}.
\end{align*}
This  leads to the same constraint for $\eta$ that we had in
\eqref{eta-bound3}.

\section*{\refname}
\bibliographystyle{plain}

\end{document}